\documentclass[11pt,a4paper]{article} 

\usepackage{amsfonts}
\usepackage{amsmath}
\usepackage{amssymb}
\usepackage{enumerate}
\usepackage{amsthm} 
\usepackage[final]{graphicx} 
\usepackage{ifthen} 
\usepackage{graphpap} 
\usepackage{longtable}
\usepackage{time} 
\usepackage{hyper} 
\usepackage{url}

\newcommand{\RootPath}{.}

\newcommand{\ExternalFiguresPath}{\RootPath/figures}

\newcommand{\eop}{\hspace*{\fill}~$\square$} 

\newtheorem{theorem}{Theorem}

\newtheorem{lemma}[theorem]{Lemma}
\newtheorem{corollary}[theorem]{Corollary}
\newtheorem{Conjecture}[theorem]{Conjecture}

\newtheorem*{theorem*}{Theorem}
\newtheorem*{proposition*}{Proposition}
\newtheorem*{lemma*}{Lemma}
\newtheorem*{corollary*}{Corollary}
\newtheorem*{Conjecture*}{Conjecture}
\newenvironment{conjecture*}{\begin{Conjecture*}}{\eop\end{Conjecture*}}

\theoremstyle{definition}
\newtheorem{Definition}[theorem]{Definition}

\newtheorem{Remark}[theorem]{Remark}
\newenvironment{remark}{\begin{Remark}}{\eop\end{Remark}}
\newtheorem{Remarks}[theorem]{Remarks} 
\newtheorem{Example}[theorem]{Example}

\newtheorem{Examples}[theorem]{Examples}

\newtheorem{Hypothesis}[theorem]{Hypothesis}

\newtheorem{Problem}[theorem]{Problem}

\newtheorem{Problems}[theorem]{Problems}

\newtheorem{Excercise}[theorem]{Excercise}

\newtheorem*{Definition*}{Definition}
\newenvironment{definition*}{\begin{Definition*}}{\eop\end{Definition*}}
\newtheorem*{Remark*}{Remark}
\newenvironment{remark*}{\begin{Remark*}}{\eop\end{Remark*}}
\newtheorem*{Remarks*}{Remarks} 
\newenvironment{remarks*}{\begin{Remarks*} {\ } \hspace{-.5cm} \begin{enumerate}}{\end{enumerate}\eop\end{Remarks*}} 
\newtheorem*{Example*}{Example}
\newenvironment{example*}{\begin{Example*}}{\eop\end{Example*}}
\newtheorem*{Examples*}{Examples}
\newenvironment{examples*}{\begin{Examples*} \begin{enumerate}}{\end{enumerate}\eop\end{Examples*}}
\newtheorem*{Hypothesis*}{Hypothesis}
\newenvironment{hypothesis*}{\begin{Hypothesis*}}{\eop\end{Hypothesis*}}
\newtheorem*{Problem*}{Problem}
\newenvironment{problem*}{\begin{Problem*}}{\eop\end{Problem*}}
\newtheorem*{Problems*}{Problems}
\newenvironment{problems*}{\begin{Problems*} {\ } \hspace{-.5cm} \begin{enumerate}}{\end{enumerate}\eop\end{Problems*}}
\newtheorem*{Exercise*}{Exercise}
\newenvironment{exercise*}{\begin{Exercise*}}{\eop\end{Exercise*}}

\newcommand{\mycaption}[1]{\centering{\vspace{\medskipamount}\refstepcounter{figure}Figure~\thefigure {#1}}}

\newcommand{\natur}{\ensuremath{\mathbb{N}}}
\newcommand{\real}{\ensuremath{\mathbb{R}}}

\newcommand{\sprod}[2]{\left<#1,#2\right>}

\newcommand{\setcond}[2]{\left\{ #1 : #2 \right\}} 

\newcommand{\equalby}[1]{\stackrel{#1}{=}}

\newenvironment{figtabular}[1]{
    \begin{figure}[htb]
    \setlength{\unitlength}{0.01\textwidth}
    \begin{center}
    \begin{tabular}{#1}
}{
    \end{tabular}
    \end{center}
    \end{figure}
}

\newcommand{\ExternalFigure}[2]{
        \includegraphics[width=#1\textwidth]{\ExternalFiguresPath/{#2}}
}

\newcommand{\conv}{\mathop{\mathrm{conv}}\nolimits}

\newcommand{\aff}{\mathop{\mathrm{aff}}\nolimits}

\newcommand{\interior}{\mathop{\mathrm{int}}\nolimits}

\newcommand{\relint}{\mathop{\mathrm{relint}}\nolimits}

\newcommand{\bd}{\mathop{\mathrm{bd}}\nolimits}

\newcommand{\cl}{\mathop{\mathrm{cl}}\nolimits}

\newcommand{\E}{\mathop{\mathbb{E}}\nolimits} 

\newcommand{\KK}{\mathop{\mathcal K}\nolimits}

\newcommand{\ortproj}{\mathbin{|}}

\newboolean{hideinwork}
\newboolean{hidefigures}
\newboolean{hideacomments}
\setboolean{hideacomments}{false} 
\setboolean{hideinwork}{false} \setboolean{hidefigures}{false}

\newcommand{\inwork}[1]{ \ifthenelse{\boolean{hideinwork}}{}{ {\bf inwork( }  #1 {\bf )} } }
\newcommand{\putfigure}[1]{ \ifthenelse{\boolean{hidefigures}}{{ \sf FIGURE }}{ #1 } }
\newcommand{\comment}[1]{}

\newcommand{\acomment}[1]{\ifthenelse{\boolean{hideacomments}}{}{ {\footnote{#1}} }} 

\setboolean{hideacomments}{true}

\renewcommand{\subset}{\subseteq} 

\newcommand{\hdmeas}{\mathcal{H}}
\newcommand{\eps}{\varepsilon}
\newcommand{\usphere}{\mathop{\mathbb{S}}\nolimits}
\newcommand{\uball}{\mathop{\mathbb{B}}\nolimits}

\newcommand{\Ccal}{\mathop{\mathcal C}\nolimits}
\newcommand{\Ucal}{\mathop{\mathcal U}\nolimits}
\newcommand{\Rot}{\mathop{\mathcal R}\nolimits_{\pi/2}}
\newcommand{\dd}{\mathrm{d}}

\begin{document}

\title{Retrieving Convex Bodies from Restricted Covariogram Functions}
\author{Gennadiy Averkov\footnote{Research supported by the \emph{Marie
Curie Research Training Network} (project \emph{Phenomena in
High Dimensions}, contract number \emph{MRTN-CT-2004-511953}) and by \emph{Deutsche Forschungsgemeinschaft} (project \emph{AV~85/1-1})} \ and Gabriele Bianchi}

\maketitle

\begin{abstract}

The covariogram $g_K(x)$ of a convex body $K \subseteq \E^d$
is the function which associates to each $x \in \E^d$ the volume of
the intersection of $K$ with $K+x.$ Matheron \cite{Matheron8601}
asked whether $g_K$ determines $K$, up to translations
and reflections in a point. Positive answers to Matheron's
question have been obtained for large classes of planar convex
bodies, while for $d\geq3$ there are both positive and negative results.

One of the purposes of this paper is to sharpen some of the known
results on Matheron's conjecture indicating how much of the
covariogram information is needed to get the uniqueness of
determination. 
We indicate some subsets of the support of the covariogram, with arbitrarily small Lebesgue measure, such that the covariogram, restricted to those subsets, identifies certain geometric properties of the body. 
These results are more precise in the planar case, but some of them, both positive and negative ones, are proved for bodies of any dimension. 
Moreover some results regard most convex bodies, in the Baire category sense.
Another purpose is to extend the class of convex
bodies for which Matheron's conjecture is confirmed by
including all planar convex bodies possessing two
non-degenerate boundary arcs being reflections of each other.

\newtheoremstyle{itsemicolon}{}{}{\mdseries\rmfamily}{}{\itshape}{:}{ }{}
\theoremstyle{itsemicolon}
\newtheorem*{msc*}{Mathematics Subject Classification (AMS 2000)}

\begin{msc*}
	52A20, 52A22, 52A38, 60D05.
\end{msc*}

\newtheorem*{keywords*}{Key words and phrases}

\begin{keywords*}
	convex body, convex polytope, covariogram, genericity result, geometric tomography, set covariance,  quasicrystal.
\end{keywords*}

\end{abstract}

\section{Introduction}

Let $K$ be a convex body in $\E^d$. The \emph{covariogram} $g_K$
of $K$ is the function
$$
g_K(x)= V(K\cap(K+x)),
$$
where $x\in\E^d$ and $V$ denotes volume in $\E^d$.
This functional, which was introduced by Matheron in his book
\cite{MR0385969} on random sets, is also sometimes called the \emph{set
covariance}, and it coincides with the  \emph{autocorrelation} of the characteristic function of $K$:
\begin{equation}\label{convoluzione}
g_K=\mathbf{1}_K\ast \mathbf{1}_{(-K)}.
\end{equation}
The covariogram $g_K$ is clearly unchanged by a translation or a reflection 
of $K$. (The term \emph{reflection} will always means reflection 
in a point.) Matheron~\cite{Matheron8601} and, independently, Adler and Pyke~\cite{MR1450931} asked the following question.

\newtheorem*{covproblem*}{Covariogram problem}

\begin{covproblem*}
Does the covariogram determine a
convex body, among all convex bodies, up to translations and
reflections?
\end{covproblem*}

Matheron conjectured a positive answer for the case $d=2,$ but this conjecture has not been completely settled.

Matheron~\cite[p.~86]{MR0385969} observed  that, 
for  $u\in \usphere^{d-1}$ and for
all $r>0$, the derivatives ${\partial}g_K(ru)/{\partial r}$ give the
distribution of the lengths of the chords of $K$ parallel to $u$.
Such information is common in \emph{stereology,} \emph{statistical shape
recognition} and \emph{image analysis,} when properties of an unknown body
are to be inferred from chord length measurements; see \cite{MR1196706},
\cite{MR1974301} and \cite{MR753649}, for example. 
Blaschke (cf. \cite{MR2162874}) asked whether the distribution of
the lengths of chords (in all directions) of a convex body characterizes the body, up to rigid motions, but Mallows and Clark~\cite{MR0259976}
proved that this is false even for convex polygons.  In fact (see
\cite{MR1232748}) the covariogram problem is equivalent 
to the problem of
determining a convex body from all its separate chord length
distributions, one for each direction $u\in \usphere^{d-1}$.

Adler and Pyke~\cite{AP91,MR1450931} asked Matheron's question in 
probabilistic terms.  Does the
distribution of the difference $X-Y$ of independent random variables
$X$ and $Y$ uniformly distributed over $K$ determine $K$, up to
translations and reflections? Since the convolution in \eqref{convoluzione} 
is, up to a multiplicative factor, the
probability density of $X-Y$, this problem is equivalent to the
covariogram one. 

Matheron's problem is also relevant in \emph{X-ray crystallography,} 
where the atomic structure of a \emph{crystal} (or \emph{quasicrystal}) is to be
found from \emph{diffraction images.}
A convenient way of describing many important examples of quasicrystals is via the \emph{cut and project scheme}. 
Here to the atomic structure, represented by a discrete set $S$ contained in a space $E$, is associated a \emph{lattice} $N$ in a higher dimensional space $E\times E'$ and a \emph{window} $W\subset E'$ (which in many cases is  a convex set). 
In this setting  $S$ coincides with the projection on $E$ of the points of the lattice $N$ which belong to $E \times W$. 
In many examples the lattice $N$ can be determined by the diffraction image. 
To determine $S$ it is however necessary to know $W$: 
the covariogram problem enters at this point, since 
the covariogram of $W$ can be obtained by the diffraction image; 
see~\cite{BaakeGrimm}.

Enns and Ehlers~\cite{MR0471018,MR967854,MR1242019} express in terms of the covariogram the distributions of random line segments in a convex body, under different types of randomness with which they are generated. The monograph \cite{MR96j:52006} contains an extensive discussion of
retrieval problems for convex bodies, while the survey
\cite{Skiena0601} deals with algorithmic aspects of reconstruction
problems in convex geometry.

The first contribution to Matheron's question was made in 1993 by
Nagel~\cite{MR1232748}, who gave a positive answer when $K$ is a planar
convex polygon; see also Schmitt~\cite{MR1196706}.
Matheron's conjecture is still
unsettled for general planar convex bodies, but it has been
confirmed  for $\Ccal^2$ convex bodies, non-strictly convex bodies, and
convex bodies that are not $\Ccal^1$; see \cite{MR2108257}.  It has been recently shown that every
convex polytope in $\E^3$ is determined by its covariogram, up to
translations and reflections (cf. \cite{Bianchi0602}). For $d \ge
4$ there exist examples of convex polytopes that are not determined
by their covariogram  (cf. \cite{MR1909913}). However \cite[p.87]{MR1416411} proves that, if $P$ is a $d$-dimensional simplicial convex polytope in general relative position with respect to $-P$, the determination by the covariogram data is unique for every $d \ge 2$ (see next section for all unexplained definitions). 
The paper
\cite{Bianchi060223}  discusses  various open retrieval problems related to the covariogram.

One of the purposes of this paper is to sharpen some of the known
results on Matheron's conjecture, indicating how much of the
covariogram information is needed to get the uniqueness of
determination.
We indicate some subsets of the support of the covariogram, with arbitrarily small Lebesgue measure, such that the covariogram, restricted to those subsets, identifies certain geometric properties of the body. 
These results are more precise in the planar case, but some of them, both positive and negative ones, are proved for bodies of any dimension. 
Moreover some results regard \emph{most} convex bodies, in the \emph{Baire category sense.}
Another  purpose is to extend the class of convex
bodies for which Matheron's conjecture is confirmed by
including all planar convex bodies possessing two
non-degenerate boundary arcs being reflections of each other.

\newcommand{\GG}{\mathbf{GC}}
\newcommand{\LC}{\mathbf{LC}}

Given two convex bodies $K$ and $H$ in $\E^d$ and a closed set $X
\subseteq \E^d,$  we introduce the following property involving $K,
H$ and $X;$ $\GG$ is a shorthand notation standing for
``covariogram coincidence'' (where covariogram is traditionally referred to by
the letter $G$).
\begin{description}
    \item[$\GG(X)$] The equality $g_K(x)=g_H(x)$ holds for all $x$ in some neighbourhood of
    $X.$
\end{description}

The following theorem
presents two choices of the set $X$ for which $\GG(X)$ implies the
coincidence of $K$ and $H$ up to translations and reflections,
under the assumption $K \in \Ccal_+^2.$ Before stating the theorem we need to introduce the notion of local symmetry and give some related explanations.

A pair of closed boundary arcs of a planar convex
body $K$ is said to be a \emph{local symmetry} of $K$ if they are reflections of each other in a point, have disjoint and nonempty relative interiors, and they are not properly contained in a pair of boundary arcs with the same properties. A planar convex body $K$ is called \emph{locally
symmetric} if it possesses a local symmetry. Planar convex bodies
without local symmetries are called \emph{globally non-symmetric.} It is known that the support of $g_K$ is the \emph{difference body} of $K$, $DK=\setcond{x-y}{x,y\in K}$, and that $DK$ is $o$-symmetric.  If $A^+$ and $A^-$ are arcs of $\bd 
K$ which compose a local symmetry, then the set $2(A^+\cup A^-)$, translated in such a way 
to be $o$-symmetric, is the union of two arcs $A$ and $-A$ of $\bd DK$. We  say that these arcs of $\bd DK$ correspond to the local symmetry $A^+, A^-$ (see Figs.~\ref{06.05.23,21:46} and \ref{05.12.19,11:36}).

A convex body $K$ is said to belong to the class $\Ccal_+^2$ if its boundary is a two-times continuously differentiable manifold and all its principal curvatures are non-zero (for detailed information see \cite[Section~2.5]{MR94d:52007}).

\newcommand{\NN}{\mathcal{N}}
\begin{theorem}\label{05.08.10,15:35}
    Let $K$ and $H$ be planar convex bodies and let $K$ be
$\Ccal_+^2$ regular. Let $\{ \pm A_n \}_{n\in \NN}$ be the collection of all the 
arcs of $\bd DK$ which correspond to local symmetries of $K$. By $x_n$ we denote the midpoint of the
segment joining the endpoints of the arc $A_n.$ 
Let $X_0:=\setcond{\pm x_n}{ n \in  \NN}$
and let $X=X_0 \cup \bd DK$ or $X=X_0 \cup \{o\}$ (see Figs.~\ref{06.05.23,21:46} and \ref{05.12.19,11:36}).
Then $\GG(X)$ implies the coincidence of $K$ and $H,$ up to
translations and reflections. \eop
\end{theorem}

        \begin{figtabular}{cc}
        \begin{picture}(40,30)
            \put(12,8){\ExternalFigure{0.14}{057/057n11v06.eps}} 
            \put(11,20){\large $K$}
        \end{picture}
        &
        \begin{picture}(40,30)
            \put(4,4){\ExternalFigure{0.30}{057/057n11v07.eps}} 
         \put(21,16){\large $o$}
         \put(3,25){\large $DK$}
        \end{picture}
        \\
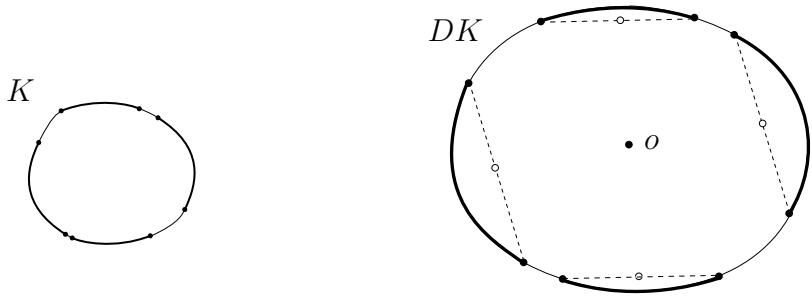
    \parbox[t]{0.40\textwidth}{\caption{The body $K$ with two local symmetries plotted in bold \label{06.05.23,21:46}}}
        &
    \parbox[t]{0.40\textwidth}{\caption{Arcs $\pm A_n$ of $DK,$ corresponding to local symmetries of $K,$ are plotted in bold; points of $X_0$ are depicted as $\circ$ \label{05.12.19,11:36}}}
    \end{figtabular}

We remark that the two choices of $X$ defined in the statement of
Theorem~\ref{05.08.10,15:35} are in some sense minimal for the assertion of the theorem to hold (for further details see Remark~\ref{06.05.16,16:17} below).
Moreover the set $X_0$ depends only on $g_K$  and not on $K$, in the sense that if $H$ and $K$ satisfy $\GG(\bd DK)$ or $\GG(o),$ then  the same set $X_0$ corresponds to $H$ and $K$. This is the content of the second part of Theorem~\ref{06.05.16,15:21}.

The following corollary is a direct consequence of
Theorem~\ref{05.08.10,15:35}.

\begin{corollary} \label{06.05.16,17:58}
    Let $K$ and $H$ be planar convex bodies and let $K$ be $\Ccal_+^2$ regular.  Then the following statements hold.
    \begin{enumerate}[I.]
    \item There exists a closed and at most countable subset $X$ of $DK,$ with no accumulation points in $\interior DK$ such that $\GG(X \cup \{o\})$ implies the coincidence of $K$ and $H,$ up to translations and reflections. Furthermore, $X$ can be chosen lying on a strictly convex curve.
    \item If $K$ is globally non-symmetric and $X$ is either $\{o\}$ or $\bd DK,$ then $\GG(X)$ implies the coincidence of $K$ and $H,$  up to translations and reflections.\eop
    \end{enumerate}
\end{corollary}

It is an open question whether Part~II of the above theorem holds for all strictly convex $K$ not necessarily $\Ccal_+^2.$ The positive answer to this question would imply the confirmation of Matheron's conjecture for all planar convex bodies.

The following theorem presents determination results involving locally symmetric and symmetric convex bodies.

\begin{theorem}\label{06.05.16,15:21} Let $K$ and $H$ be planar convex bodies with $K$ strictly convex. Then the following statements hold true.
    \begin{enumerate}[I.]
        \item If $K$ is $o$-symmetric and $A$ is a closed simple curve in $\interior DK$ bounding an open region in $DK$ that contains the origin, then the conditions $DK=DH$ and $\GG(A)$ imply the coincidence of $K$ and $H,$ up to translations and reflections.
        \item If $X$ is either $\{o\}$ or $\bd DK,$ then $\GG(X)$ implies the coincidence of all local symmetries of $K$ and $H,$ up to translations.     \eop
    \end{enumerate}
\end{theorem}

The statement of the following corollary follows directly from the
second part of Theorem~\ref{06.05.16,15:21} and two results from
\cite[Theorem~1.1 and Proposition~1.4]{MR2108257}.

\begin{corollary} \label{06.05.24,09:59}
    Every locally symmetric convex body in $\E^2$ is determined uniquely, up to translations and reflections, by its covariogram function. \eop
\end{corollary}

Let $\uball^d$ denote the closed unit ball in $\E^d$ centered at the origin. Further on, we introduce the condition $\GG'(X),$ a relaxation of $\GG(X),$ and the ``local coincidence'' condition $\LC.$

\begin{description}
    \item[$\GG'(X)$] The equality $g_K(x)=g_H(x)+c$ holds for all $x$ in some neighbourhood of
    $X$ and a suitable constant $c \in \real.$
    \item[$\LC$] For every boundary point $p$ of $K$ there
    exists a boundary point $q$ of $H$ such that for some $\eps
    >0$ the bodies $K \cap (p+\eps \uball^d)$ and $H \cap (q+ \eps
    \uball^d)$ coincide, up to translations and reflections; the
    same statement also holds with the roles of $K$ and $H$
    interchanged.
\end{description}

In the following theorem the relationship between the conditions $\GG(\{o\}),$ $\GG'(\{o\}),$ $\GG(\bd DK),$ and $\LC$ is discussed.

\begin{theorem}\label{06.05.15,12:36}
    Let $K$ and $H$ be convex bodies in $\E^d, \ d \ge 2.$ Then the following statements hold true.
    \begin{enumerate}[I.]
        \item If $K$ is strictly convex, then $\LC$ implies $\GG(\bd DK).$
        \item If $d=2$ and $K$ is $\Ccal_+^2,$ then $\LC$ is equivalent to $\GG(\bd DK).$
        \item If $d=2$ and $K$ is strictly convex, then $\GG'(\{o\})$ is equivalent to $\GG(\bd DK).$
        \item There exist planar convex bodies $K$ and $H$ belonging to the class $\Ccal_+^2$ such that  $\GG(\bd DK)$ holds, while  $\GG(\{o\})$ does not.
        \item There exist  $K$ and $H,$ which are convex $d$-polytopes,  such that the conditions $DK=DH$ and $\GG(\{o\})$ hold, while $\GG(\bd DK)$ does not.     \eop
    \end{enumerate}
\end{theorem}

It is an open problem whether for $\Ccal_+^2$ convex bodies in
$\E^d,$ $d \ge 3,$ the condition $\GG(\bd DK)$ implies $\LC.$

 The space $\KK^d$ endowed with the Hausdorff metric is
locally compact and by this a Baire space (see \cite{MR1243011}
and \cite[p.119]{MR94d:52007}). Thus, we may speak about 
statements that hold for most convex bodies, i.e., for all convex
bodies with at most a meager set of exceptions. We recall that a set
is said to be \emph{meager} if it is a countable union of nowhere
dense sets and \emph{residual} if it is a complement of a meager
set. Trivially, a finite intersection of residuals is a again a
residual. Furthermore, every set possessing a residual subset is
also a residual.

\begin{theorem}\label{06.05.08,17:28}
    In the Baire category sense, for most convex bodies $K$ and all convex bodies $H$ in $\E^d$ the following statements hold true.
    \begin{enumerate}[I.]
        \item If $K$ and $H$ satisfy $\GG(\bd DK),$ then $K$ and $H$ coincide, up to translations and reflections.
        \item If $d=2$ and $K$ and $H$ satisfy $\GG(\{o\}),$ then $K$ and $H$ coincide, up to translations and reflections.
    \end{enumerate}
    Furthermore, the above two statements cannot be extended to all pairs of convex bodies $K$ and $H,$ since there exist bodies $K$ and $H$ not coinciding, up to translations and reflections, and such that $\GG(\bd DK \cup \{o\})$ holds. In addition, the bodies $K$ and $H$ satisfying the above conditions can be chosen belonging to the class $\Ccal_+^2.$ \eop
\end{theorem}

It is an open question whether  Part~II of the statement can be carried over to convex bodies of higher dimensions.

The previous theorem is related to \cite[Theorem~2]{MR1416411} and  \cite[Theorem~6.2]{MR1938112}. 
Theorem~2 from \cite{MR1416411} states that most convex bodies $K$ in $\E^d$, for any $d\geq3$, are determined by the combined knowledge of the \emph{width} of $K$ in direction $u$ and of the $d-1$-dimensional volume of $K \ortproj u^\perp$, for each $u\in \usphere^{d-1}$. 
Here $u^\perp$ denotes the orthogonal complement to $u,$ and $K \ortproj u^\perp$ stands for the orthogonal projection of $K$ onto the hyperplane $u^\perp.$ (See also \cite{MR1610155}, where
 ``local versions'' of some results from \cite{MR1416411} are
 obtained,  as well as  \cite[Theorems~3.3.17 and
 3.3.18]{MR96j:52006}.)
The relation with the covariogram comes from the fact that knowing its support $D K$ is equivalent to knowing the width of $K$ in all directions, and that the knowledge of $g_K$ in a neighbourhood of $o$ gives the volumes of all $(d-1)$-dimensional projections of $K$. This follows from the formula
\begin{equation} \label{05.11.07,11:34}
    \left. \frac{\partial^+ g_K(ru)}{\partial r} \right|_{r=0} = -V_{d-1}(K \ortproj u^\perp)\quad (u\in \usphere^{d-1}),
\end{equation}
 proved in \cite{Matheron8601} and \cite{MR0385969}. Here $\partial^+ / \partial r$ stands for right derivative.  
Theorem~6.2 from \cite{MR1938112} is another result related to Theorem~\ref{06.05.08,17:28}, which states that most convex planar bodies are
 determined by the covariogram function over its entire support. 

In Section~3 we prove Theorem~\ref{06.05.15,12:36}. Sections~4 and 5 are independent of each other and present proofs of
Theorem~\ref{05.08.10,15:35} and Theorem~\ref{06.05.16,15:21}, respectively. In Section~6 we prove Theorem~\ref{06.05.08,17:28}.

\section{Background from convex geometry}

The Euclidean $d$-dimensional space with the origin $o,$ scalar product $\sprod{\,.\,}{\,.\,},$ and the
norm $|\, . \,|$ is denoted by $\E^d.$ The unit sphere and the
unit ball in $\E^d$ are denoted by $\usphere^{d-1}$ and
$\uball^d,$ respectively. The orthogonal projection of a set $X
\subseteq \E^d$ onto an affine space $L \subseteq \E^d$ is denoted
by $X | L.$ If $u \in \E^d,$ then $u^\perp$ stands for the
orthogonal complement to $u.$ By $V_j, \ j \in \{1,\ldots,d\},$ we denote the $j$-dimensional volume of a convex set in $\E^d$ of dimension at most $j.$ If $j=d,$ we omit
the subscript and write simply $V.$ We write $\hdmeas^j$ for the $j$-dimensional Hausdorff measure in $\E^d.$ In what follows, in integrals on $j$-dimensional spheres in $\E^d$ (with $j \le d-1$) we abbreviate $\hdmeas^j(\dd u)$ by $\dd u.$ 

The abbreviations $\bd,$
$\interior,$ $\relint,$ $\cl$ and $\aff$ stand for boundary, interior, relative interior, closure and affine hull,
respectively.

Following the
monograph \cite{MR94d:52007} we denote by $\KK^d$ and $\KK^d_0$
the classes of non-empty, compact, convex sets and compact, convex
sets with non-empty interior, respectively. Elements of $\KK^d_0$ are said to be \emph{convex bodies.}  If $o \in K,$ then
$r_K(u):=\max \setcond{\alpha>0}{\alpha u \in K}, \ u \in
\usphere^{d-1},$ is called the \emph{radius function} of $K.$ We also introduce the \emph{support function} of $K$ by $h_K(u):=\max \setcond{\sprod{x}{u}}{x \in K}, \ u \in \E^d.$ The \emph{difference body} of
a convex body $K$ is the set $DK:=K+(-K).$ The function $w_{K}(u):=h_{DK}(u)$ is called the \emph{width function} of $K.$ Observe that for $u \in \usphere^{d-1}$ the quantity $w_K(u)$ is equal to the distance between the two distinct supporting hyperplanes of $K$ orthogonal to $u.$  The face of a convex body $K
\subseteq \E^d$ in the direction $u \in \E^d \setminus \{o\}$ is
denoted by $F_K(u).$ Elements of $\KK^d_0$ representable as intersection of a finite number of closed halfspaces are called \emph{convex $d$-polytopes.}
Two convex $d$-polytopes $P, Q \subseteq \E^d$ are said to be in a \emph{general relative position} if for every $u \in \usphere^{d-1}$ and $x \in \E^d,$ $F_P(u) \cap (F_Q(u)+x)$ is either empty or a singleton.

If $K$ is a convex body in $\E^d$ and $u \in \E^d \setminus \{o\},$ then,  for $p \in \relint F_K(u),$ the normal cone $N_K(p)$ does not depend on the choice of $p$ in $\relint F_K(u).$ We denote $N_K(p)$ by $\overline{N}_K(u)$ and call it the \emph{normal cone of $K$ in direction $u.$}

 If $X$ is a subset of
$\bd K,$ then the set of all outward unit normals of $K$ at points
of the set $X$ is called the \emph{spherical image} of $X$
with respect to $K.$ Two boundary arcs $A$ and $B$ of $K \subseteq
\E^2$ are said to be \emph{antipodal} if their spherical
images with respect to $K$ are reflections of each other.  

For $\eps>0$ the \emph{$\varepsilon$-neighbourhood of a set} $X \subseteq \E^d$ is the
set $X+ \varepsilon \cdot \interior \uball^d,$ that is, the set
consisting of all those points $x$ whose distance to some point of
$X$ is strictly less than $\varepsilon.$  The Hausdorff distance
$\delta(X,Y)$ between non-empty compact sets $K$ and $H$ in $\E^d$
is the least possible $\alpha\ge 0$ such that $X \subseteq Y +
\alpha \cdot \uball^d$ and $Y \subseteq X + \alpha \cdot
\uball^d.$ Information on the Hausdorff distance in the class of
convex bodies is collected in \cite[\S1.8]{MR94d:52007}. We introduce the  distance function $\bar{\delta}(X,Y)$ for sets $X,Y \subseteq \E^d$ as the minimum of $\delta(X,\phi(Y)),$  where $\phi$ ranges over all translations and reflections.

The \emph{area measure of order $d-1$} of a convex body $K\subseteq \E^d$ (see \cite[\S4.2]{MR94d:52007}) is given by 
\begin{equation*}
   S_{d-1}(K,\omega):=\hdmeas^{d-1}(\setcond{p \in \bd K}{\mbox{some outward normal of $K$ at $p$ belongs to $\omega$}}),
\end{equation*}
 where $\omega$ is a Borel set in $\usphere^{d-1}.$ If $d=2,$ then $S_1(K,\,.\,)$ is said to be the \emph{length measure} of $K.$ 

Given a strictly convex body $K \subseteq \E^d,$  $z_K(u), \ u \in
\usphere^1,$ denotes the boundary point of $K$ with outward
normal $u.$ 

If $S_{1}(K,\,.\,)$ is absolutely continuous with respect to the
measure $\hdmeas^{1}$ on $\usphere^{1},$ then we denote by $R_K$ the
\emph{Radon-Nikodym derivative} of $S_{1}(K,\,.\,),$ i.e.\ 
the function obeying the equality
\begin{equation}
    S_{1}(K,\omega)=\int_{\omega} R_K(u) \, \hdmeas^{1}(\dd u)
\end{equation}
for all Borel subsets $\omega$ of $\usphere^1$ (see also Formula~(4.2.20) from \cite{MR94d:52007}).  If  the function $R_K$ is continuous, $R_K(u), \ u \in \usphere^1,$ is
the radius of curvature of $K$ at the point $z_K(u).$
We remark that, according to our definition, $R_K(u)$ may be zero
for some values of $u.$ If $K \in \Ccal_+^2,$ then $R_K$ is continuous and strictly positive, and
$\tau_K(u)=1/R_K(u)$ is said to be the \emph{curvature} of $K$ at $z_K(u).$

We parametrize the unit circle $\usphere^1$ in a standard manner
by the vector function $u(t) :=(\cos t, \sin t),$ where $t \in
\real.$  Given a $\Ccal_+^2$ planar convex body $K,$ we
put $\tau_K(t):=\tau_K(u(t)),$ $z_K(t):=z_K(u(t))$ and $R_K(t):=R_K(u(t)).$ For $t_1,
t_2 \in \real$ with $t_1 \le t_2$ we introduce the notation
$z_K(t_1,t_2):= \setcond{ z_K(t)}{ t_1 \le t \le t_2}.$ For
$\Ccal_+^2$ convex bodies $K, H \subseteq \E^2$ one can show that
$z_{K+H}(t)=z_K(t)+z_H(t).$  It is known that for $\Ccal_+^2$ planar convex
bodies $K$ and $H$ one also has $R_{K+H}(t)=R_K(t)+R_H(t).$

The knowledge of $z_K(t_1)$ and the curvature $\tau_K(t)$ for $t
\in [t_1,t_2]$ allows to determine the arc $z_K(t_1,t_2).$ More
formally, the parametrization $z_K(t), \ t_1 \le t \le t_2,$ of
$z_K(t_1,t_2)$ is determined from the representation
\begin{equation} \label{05.02.22,11:26}
    z_K(t) = z_K(t_1) + \int_{t_1}^{t} \frac{u(s+\pi/2)}{\tau_K(s)} \, \dd s,
\end{equation}
Equality \eqref{05.02.22,11:26} can be
found in \cite[p.~11]{MR840401} and \cite[p.~186]{MR1938112} (cf.\  also a more general result
stated in Theorem~4.3.2 from \cite{MR94d:52007}). By
\eqref{05.02.22,11:26} we see that two antipodal arcs
$z_K(t_1,t_2)$ and $z_K(t_1+\pi,t_2+\pi)$ with $t_1<t_2$ and
$t_2-t_1<\pi$ are reflections of each other if and only if
$\tau_K(t)=\tau_K(t+\pi)$ for every $t \in [t_1,t_2].$

\section{Relationship between $\LC$ and covariogram data}

In this section we prove Parts~I, III, IV and V of Theorem~\ref{06.05.15,12:36}. The proof of Part~II is postponed to the next section.

\begin{proof}[Proof of Theorem~\ref{06.05.15,12:36} (Part~I)]
Using $\LC$ and strict convexity of $K,$ we see that for every $u \in \usphere^{d-1},$ there exists a relatively open subset $G_u$ of $\usphere^{d-1}$ such that $(-G_u) \cap G_u = \emptyset,$ $u \in G_u,$ and the equality 
\begin{equation} \label{06.09.21,17:00}
	\{z_K(G_u)+x_1, -z_K(-G_u)-x_2\} = \{z_H(G_u), -z_H(-G_u)\}
\end{equation}
holds for some $x_1, x_2 \in \E^d$ depending on $u.$ Then, since $z_{DL}(v)=z_{L}(v)-z_L(-v)$ for every $L \in \KK_0^d$ and $v \in \usphere^{d-1},$  from \eqref{06.09.21,17:00} we deduce
\begin{equation}\label{some_label}
z_{DK}(G_u)+x_1-x_2=z_{DH}(G_u).
\end{equation}
Since $\usphere^{d-1}$ is compact, there is a finite sub-family $\setcond{\pm G_i}{i=1,\ldots,n} , \ n \in \natur,$ of $\setcond{\pm G_u}{u \in \usphere^{d-1}}$ covering $\usphere^{d-1}.$ By \eqref{some_label} we have
\begin{equation} \label{06.09.11,12:18}
	S_{d-1}(DK,\omega)=S_{d-1}(DH,\omega)
\end{equation} 
for every Borel set $\omega$ being a subset of $G_i$ or $-G_i$ for some $i=1,\ldots,n.$  
Furthermore, \eqref{06.09.11,12:18} can be derived for any Borel subset $\omega$ of $\usphere^{d-1}$ using the decomposition 
$$
	\omega = \bigcup_{i=1}^n (\omega \cap G_i) \cup (\omega \cap (-G_i))
$$
and the inclusion-exclusion principle. Hence $S_{d-1}(DK,\, . \,)=S_{d-1}(DH,\, . \,),$ and by this $DK=DH$ (and $x_1=x_2$ for each $u \in \usphere^{d-1}$). Using the latter equality together with \eqref{06.09.21,17:00}, we deduce that for every $u \in \usphere^{d-1}$ the functions $g_K$ and $g_H$ coincide in some neighbourhood $X(u)$ of $z_{DK}(u).$ Consequently, $g_K$ and $g_H$ coincide in the open set $\bigcup_{u \in \usphere^{d-1}} X(u),$ which encloses $\bd DK.$
\end{proof}

For $x\in \interior DK \setminus \{o\}$ let $p$ and $q$ be the
endpoints of the arc $(K+x)\cap \bd K$. Then $p-x$ and $q-x$ belong to $\bd
K$ and $P_K(x):=\conv \{ p,q,p-x,q-x\}$ is a parallelogram. Following
\cite{Matheron8601} we define $D_K(x)$ to be
$\pm (p-q)$ with the sign determined by the condition $\sprod{x}{\Rot
D_K(x)}<0,$ see Figs.~\ref{05.12.30,01:24} and \ref{06.09.22,09:56}. It is known that
\begin{equation} \label{06.09.21,17:07}
	D_K(x)=\Rot(\nabla g_K(x)),
\end{equation}
where $\Rot$ denotes, throughout the paper, the counterclockwise rotation about the origin by the angle $\pi/2.$ 
 We also have the equality
\begin{equation}\label{iteration}
D_K(D_K(x))=-x.
\end{equation}

    \begin{figtabular}{cc}
        \begin{picture}(34,25)
        \put(-1,-0.7){\ExternalFigure{0.33}{057/057n01v02.eps}}
        \put(32,19){$K$}
        \put(26,1){$p$}
        \put(26,22){$q$}
        \put(-1,-2){$p-x$}
        \put(-1.5,18){$q-x$}
        \end{picture}
      &
        \begin{picture}(34,25)
         \put(3,0){\ExternalFigure{0.33}{057/057n01b.eps}}
	\put(31,6){$x$}
	\put(9,17){$D_K(x)$}
	\put(6,-0.5){$o$}
        \end{picture}	
     \\
    \parbox[t]{0.45\textwidth}{\mycaption{\label{05.12.30,01:24}}}
     &
    \parbox[t]{0.45\textwidth}{\mycaption{\label{06.09.22,09:56}}}
    \end{figtabular}

\begin{proof}[Proof of Theorem~\ref{06.05.15,12:36} (Part~III)]
First we prove that $\GG'(\{o\})$ implies $\GG(\bd DK).$ In the plane $V_1(K \ortproj u^\perp)=w_{K}(\Rot u) = h_{DK}(u)$ for every $u \in \usphere^1.$ Hence, from \eqref{05.11.07,11:34} we get that  $\GG'(\{o\})$ implies $DK=DH.$

It is easy to prove that the mapping $D_K : \interior DK \setminus \{o\} \rightarrow \interior DK \setminus \{o\}$ is continuous and, by
\eqref{iteration}, also its inverse is continuous. Moreover it
maps a punctured neighbourhood $G$ of $o$ in a set  $G'$ with $\bd
DK \subseteq \bd G'$. Let $G$ be a punctured neighbourhood of $o$
in which $g_K$ and $g_H$ coincide, up to an additive constant.
Then, in view of \eqref{06.09.21,17:07},  $D_K(x)=\Rot(\nabla g_K(x))=\Rot(\nabla g_H(x)) = D_H(x)$ for all $x \in G$
and we get that $P_K(x)$ and $P_H(x)$ are translates of each other
(cf.\  the related Lemma 1.5 from \cite{MR2108257}).

The set $K \setminus \interior P_K(x)$ is the  union of the four lunettes outside
$P$ (where a \emph{lunette} of $K$ is a compact set bounded by a
chord of $K$ and a boundary arc of $K$ joining the endpoints of
this chord). Let $p$ and $q$ be as in the definition of $D_K.$ The sum of the areas of the two lunettes adjacent to
$[p,q]$ and $[p,q]-x$ equals $g_K(x)$, the sum of the areas of the
other two lunettes equals $g_K(D_K(x))$. Similar considerations hold
for $H,$ too, and thus
\begin{eqnarray}\label{sums}
g_K(D_K(x)) &=& V(K)-V(P_K(x))-g_K(x)=g_K(o)-V(P_K(x))-g_K(x), \\
g_H(D_H(x))&=& V(H)-V(P_H(x))-g_H(x)=g_H(o)-V(P_H(x))-g_H(x). \label{sums2}
\end{eqnarray}
Since $g_K(x)-g_K(o)=g_H(x)-g_H(o),$ for $x \in G,$ we obtain from
\eqref{sums} and \eqref{sums2} that  $g_K(x)=g_H(x)$ for  $x \in G'.$

The proof of the converse implication is similar. We argue
backwards: an open subset $G'$ of $DK$ with $\bd G' \subseteq \bd
DK$ and $g_K(x)=g_H(x),$ for $x \in G,'$ is mapped by $D_K$ onto a
punctured neighbourhood $G$ of the origin. Thus, using
\eqref{sums} and \eqref{sums2}, we obtain $g_K(x)-g_H(x)=V(K)-V(H)$ for $x \in G.$
\end{proof}

\begin{proof}[Proof of Theorem~\ref{06.05.15,12:36} (Part~IV)]
In view of Theorem~\ref{06.05.15,12:36} (Parts~I and II) it suffices to construct $K$ and $H$ satisfying $\LC$ and $V(K)<V(H).$ Let $T$ be a regular triangle of unitary edge length with center
at the origin. If $u_1,u_2,u_3$ denote the unit outer normals to the edges, the area measure of $T$ is $\sum_{i=1}^3 \delta_{u_i} $, where $\delta_{u_i}$ is the Dirac delta distribution on the manifold $\usphere^1$ centered in $u_i$.
For each $i=1,2,3$ let $\phi_i$ be a continuous non-negative function on $\usphere^1$ supported in a small arc centered at $u_i$,  whose integral on $\usphere^1$ is  $1$.
Moreover choose $\phi_i$'s in such a way  that $\int_{\usphere^1}\sum_{i=1}^3  \phi_i(u)u \, \dd u=0.$ Let $K_1$ be a convex body in $\E^2$ whose length measure has density $\sum_{i=1}^3\phi_i.$
It is clear that the measure $S_1(K_1,\, . \,)$ approximates in some sense the measure $S_1(T,\, . \,).$ In fact, it can be shown that the  \emph{Prohorov metric} (see \cite{MR1962583} for the definition) of $S_1(K_1,\, . \,)$ and $S_1(T,\, . \,)$ can be made arbitrarily small. 

Let $K_2$ be a
slight rotation of $K_1$ such that the supports of the area
measures of $K_1$ and $K_2$ are disjoint. By a stability result for the Minkowski problem with respect to the Prohorov metric proved in \cite[Theorem~3.1]{MR1962583}  $K_1+K_2$ is close to $2T$, while $K_1-K_2$ is close to $DT$ in the metric $\bar{\delta}$ (Theorem~3.1 from \cite{MR1962583} is a strengthening of Theorem~7.2.2 from \cite{MR94d:52007}, see also related Theorems~4.3.5 and 7.2.6 from this monograph). Consequently, the area of
$K_1+K_2$ is close to $4 V(T),$ while the area of $K_1-K_2$ is
close to $V(DT)=6 V(T)$. Let
    \begin{align*}
       &K :=K_1 + K_2 + \uball^2,& 
       &H :=K_1 - K_2 + \uball^2.
    \end{align*}
See Figs.~\ref{06.01.04,09:52} and \ref{06.01.04,09:53} depicting possible choices of bodies $K$ and $H,$ respectively.

By construction, $K$ and $H$ are $\Ccal_+^2$ and they satisfy $\LC.$
Thus, by Theorem~\ref{06.05.15,12:36} (Parts~I and II), they have equal
covariograms in a neighbourhood of $\bd DK.$ Furthermore, $K$ has
smaller area then $H,$ since by the two-dimensional version of the \emph{Steiner formula} (see \cite[Section~4.1]{MR94d:52007}) we have $V(K)=V(K_1+K_2)+V_1(\bd K_1)+V_1(\bd
K_2)+\pi$ and $V(H)=V(K_1-K_2)+V_1(\bd K_1) + V_1(\bd K_2)+\pi.$
\end{proof}
        \begin{figtabular}{cc}
        \begin{picture}(30,25)
            \put(-20,-20){\ExternalFigure{0.60}{057/057n16.eps}} 
        \put(3,20){$K$}
        \end{picture}
        &
        \begin{picture}(30,25)
            \put(-20,-30){\ExternalFigure{0.60}{057/057n17.eps}} 
        \put(3,20){$H$}
        \end{picture}
        \\
    \parbox[t]{0.45\textwidth}{\mycaption{\label{06.01.04,09:52}}}
        &
    \parbox[t]{0.45\textwidth}{\mycaption{\label{06.01.04,09:53}}}
    \end{figtabular}

In the proof of Theorem~\ref{06.05.15,12:36} (Part~V) we shall need the
following lemma, presenting a formula which is also related
to a formula given in \cite[p.18]{nagel:habil}.

\begin{lemma} \label{06.01.06,01:19}
    Let $P$ be a convex polygon in $\E^2,$ and $G(P)$ be given by
    \begin{equation*}
        G(P):= \bigcap \setcond{DT}{T=\conv \{p_1,p_2,p_3\}, \ p_1, p_2,
        p_3 \ \mbox{are consecutive vertices of} \ P}.
    \end{equation*}
    Then $o \in \interior G(P),$ and  for $u \in G(P)$ we have
    \begin{equation*}
        g_P(u) = V(P) - w_P(\Rot u) + |u|^2 C,
    \end{equation*}
    where  $C$ depends only on
$$\{\overline{N}_P(\Rot u), \overline{N}_P(-\Rot u)\}.$$
\end{lemma}
\begin{proof}
Let $p_1$ and $p_2$ be antipodal vertices of $P$ such that $u \in N_P(p_1)$ and $-u \in N_P(p_2).$  Let $I_j$ be the union of the two edges of $P$ adjacent to $p_j \ (j=1,2).$
    In view of the assumption $u\in G(P)$ we have  $I_j \cap (I_j+u) \ne \emptyset$. The closure of the set
 $(P + [o,u]) \setminus (P \cup(P+u))$ consists of two triangles $\Delta_j, \ j=1,2$, possibly empty,  where one edge of $\Delta_j$ is $[p_j,p_j+u]$ and the other two edges are  parallel to the two edges of $I_j$ (see Fig.~\ref{06.01.06,11:10}). Therefore  $V(\Delta_1)+V(\Delta_2)=C |u|^2$
where $C$ depends only on the directions of the edges of $P$ adjacent to $p_1$ and to $p_2.$
Since
    \begin{equation*}
        g_P(u)= V(P \cap (P+u)) = V(P+u) - V((P+[o,u])\setminus P) +
        V(\Delta_1) + V(\Delta_2)
    \end{equation*}
    and $V((P+u)\setminus P)= w_P(\Rot u)$ we get the desired formula.
\end{proof}

        \begin{figtabular}{c}
        \begin{picture}(40,44)
            \put(4,4){\ExternalFigure{0.30}{057/057n19.eps}} 
        \put(7,30){$P$}
        \put(30,30){$P+u$}
        \put(17,5){$\Delta_2$}
        \put(16,33){$\Delta_1$}
	\put(13.8,16.7){$P \cap (P+u)$}
        \end{picture}
        \\
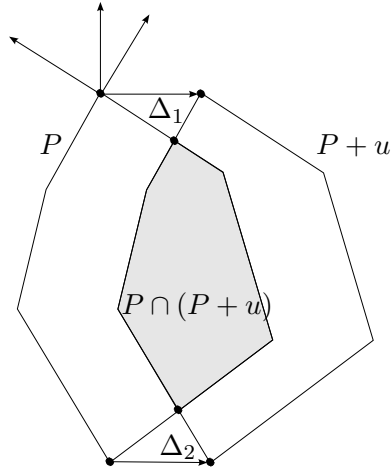
    \parbox[t]{0.45\textwidth}{\caption{$g_K(u)$ is the area of the filled region\label{06.01.06,11:10}}}
    \end{figtabular}

    If $d_i \in \natur$ and $K_i$ is a convex body in $\E^{d_i} \ (i=1,2),$ then the covariogram function of $K:=K_1 \times K_2$ is given by
    \begin{equation} \label{05.12.27,12:24}
        g_K(x)=g_{K_1}(x_1)g_{K_2}(x_2),
    \end{equation}
 where $x_i \in \E^{d_i} \ (i=1,2)$ and $x:=(x_1,x_2).$

\begin{proof}[Proof of Theorem~\ref{06.05.15,12:36} (Part~V)]
  We introduce convex polygons $P_1, P_2 \subseteq \E^2$ which are
  obtained from the square $Q:=[-10,10]^2$ by ``cutting off''
  isosceles triangles at the vertices of $Q.$ The polygon $P_k \
  (k=1,2)$ is constructed by cutting off the isosceles triangle
  with lateral sides having length $\alpha_{i,j}^k$ at the vertex $i e_1 + j e_2$
  of $Q$ for each $i, j \in \{-1,1\},$ where the constants
  $\alpha_{i,j}^k$ are defined as follows:
  \begin{align*}
    \alpha_{1,1}^1 &=10 & \alpha_{1,1}^2 &=9 \\
    \alpha_{-1,1}^1 &=2 & \alpha_{-1,1}^2 &=1 \\
    \alpha_{-1,-1}^1 &=2 & \alpha_{-1,-1}^2 &=3 \\
    \alpha_{1,-1}^1 &=8 & \alpha_{1,-1}^2 &=9
  \end{align*}
  See Figs.~\ref{05.12.27,11:43} and \ref{05.12.27,11:44}
  depicting $P_1$ and $P_2.$
        \begin{figtabular}{cc}
        \begin{picture}(24,25)
            \put(2,2){\ExternalFigure{0.20}{057/057n12.eps}} 
        \put(1,23){$P_1$}
        \end{picture}
        &
        \begin{picture}(24,25)
            \put(2,2){\ExternalFigure{0.20}{057/057n13.eps}} 
        \put(1,23){$P_2$}
        \end{picture}
        \\
    \parbox[t]{0.45\textwidth}{\mycaption{\label{05.12.27,11:43}}}
        &
    \parbox[t]{0.45\textwidth}{\mycaption{\label{05.12.27,11:44}}}
    \end{figtabular}
No translation or reflection of $P_1$ coincides with $P_2,$ and moreover, for each $k \in \{1,2\}$,
\begin{equation*}
    V(P_k) = V(Q)- \frac{1}{2} \sum_{i,j=-1}^1 (\alpha_{i,j}^k)^2 = 20^2 -
    \frac{1}{2} \cdot 172
\end{equation*}
Furthermore, it is easy to see that $DP_1=DP_2.$
Lemma~3.1 from
\cite{MR1909913} proves that  for each $u\in \usphere^1$ the knowledge of the covariogram of a convex polygon $P \subseteq \E^2$ near the boundary of its support determines the set  $\{V_1(F_P(-u)), V_1(F_P(u))\}.$
But for $u=(1,1)$ we have 
$$
	\{10 \sqrt{2}, 2 \sqrt{2} \} = \{V_1(F_{P_1}(-u)), V_1(F_{P_1}(u))\} \ne \{V_1(F_{P_2}(-u)), V_1(F_{P_2}(u))\} = \{9 \sqrt{2}, 3 \sqrt{2}\}.
$$
Hence $g_{P_1}$ and $g_{P_2}$ do not coincide in some neighbourhood of $\bd DP_1=\bd DP_2.$

For $j \in \{1,2\}$, let $G(P_j)$ be   defined as in the statement of Lemma~\ref{06.01.06,01:19}.  Since for each $u\in \real^2$
$$
\{\overline{N}_{P_1}(\Rot u), \overline{N}_{P_1}(-\Rot u)\}=\{\overline{N}_{P_2}(\Rot u), \overline{N}_{P_2}(-\Rot u)\},
$$
 Lemma~\ref{06.01.06,01:19} implies $g_{P_1}(v)=g_{P_2}(v)$ for each $v\in G(P_1)\cap G(P_2)$. Since this set is a  neighbourhood of $o$ the proof for $d=2$ is concluded by putting $K=P_1$ and $H=P_2.$

For $d \ge 3$ we define $K =
 P_1 \times [-1,1]^{d-2}, \ H = P_2 \times [-1,1]^{d-2},$ and  the property
 \eqref{05.12.27,12:24} proves the assertion.
\end{proof}

\section{Determination results for planar $\Ccal_+^2$ bodies}

In \cite[Lemma~2.1]{MR1938112} it is shown that for planar
$\Ccal_+^2$ convex bodies the asymptotic behaviour of the covariogram
function near the boundary of its support allows to determine the
non-ordered pair $\{\tau_K(u),\tau_K(-u)\}$ for every $u \in
\usphere^1.$ Thus, the following lemma holds.

\begin{lemma} \label{05.08.10,10:47}
    Let $K$ be a planar $\Ccal_+^2$ convex body. Then the
covariogram of $K$ over any neighbourhood of $\bd DK$ determines
the mapping $u \mapsto \{\tau_K(u),\tau_K(-u)\},$ where $u$ ranges
over $\usphere^1.$ \eop
\end{lemma}

Suppose that $A$ and $B$ are two disjoint antipodal boundary arcs of
$K.$ Let $z$ be an endpoint of $B.$ Let us denote by $\bar{B}$ the
convex curve obtained by joining $B$ and  the appropriate half of
the line which is tangent to $B$ at $z.$ We say that the
translated arc $A+u, \ u \in \E^d,$ \emph{captures} the endpoint $z$ of $B$ if $A$
intersects $\bar{B}$ at two points which bound an arc of $\bar{B}$
containing $z$ in its relative interior (see Fig.~\ref{05.12.16,17:43}).

        \begin{figtabular}{cc}
        \begin{picture}(28,20)
            \put(4,4){\ExternalFigure{0.20}{057/057n09.eps}} 
            \put(6,17){$A$} %
            \put(13,2){$A+u$} %
            \put(18,6){$z$} %
            \put(25,12){$B$} %
        \end{picture}
	&
        \begin{picture}(30,24)
            \put(0,-2){\ExternalFigure{0.30}{057/057n36.eps}} 
	   \put(12,-2){\small $u(t_1)$}
	   \put(21,22){\small $u(t_1+\pi)$}
	   \put(29,8){\small $(\alpha',\beta)$}
        \end{picture}
        \\
    \parbox[t]{0.45\textwidth}{\mycaption{\label{05.12.16,17:43}}}
   &
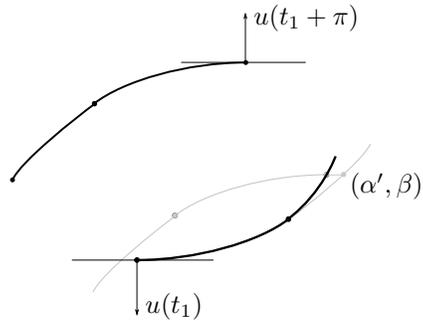
    \parbox[t]{0.45\textwidth}{\caption{$z_K(t_1,t_2)$ and $z_K(t_1+\pi,t_2+\pi)$ are bold lines; translation and reflection of $z_K(t_1+\pi,t_2+\pi)$ are gray lines\label{06.09.22,10:43}}}
    \end{figtabular}

The following lemma on capturing arcs improves slightly Lemma~4.2 from
\cite{MR1938112}, since it also indicates which translation vector
can be chosen for making a capture.
\begin{lemma} \label{05.08.10,14:15}
    Let $K$ be a planar ${\Ccal}_+^2$ convex body.   
Assume that the antipodal arcs  $z_K(t_1,t_2)$ and
$z_K(t_1+\pi,t_2+\pi)$ (where $t_1<t_2$ and $t_2-t_1<\pi$) are not
reflections of each other. Let $t^\ast$ be equal to $t_1$ for the
case $\tau_K(t_1)\ne\tau_K(t_1+\pi)$ and be equal to the maximal value  in
$[t_1,t_2]$ such that for every $t \in [t_1,t^\ast]$ the equality
$\tau_K(t)=\tau_K(t+\pi)$ holds, otherwise. We put
$z_0:=\frac{1}{2} ( z_{DK}(t_1)+z_{DK}(t^\ast)).$  Then there
exists a vector $u$ arbitrarily close to $z_0,$ such that either
$z_K(t_1)$ is captured by  $z_K(t_1+\pi,t_2+\pi)+u$  or
$z_K(t_1+\pi)$ is  captured by $z_K(t_1,t_2)-u.$
\end{lemma}
\begin{proof}
Without loss of generality we assume that $u(t_1)=(0,-1)$ and
$z_K(t_1)=o.$  The mapping $R z := -z + z_K(t_1)+z_K(t_1+\pi)$ is
a reflection with respect to the midpoint of the segment
$[z_K(t_1),z_K(t_1+\pi)].$ By the choice of $t^\ast$ and by
\eqref{05.02.22,11:26} we see that the arcs $z_K(t_1,t)$ and
$Rz_K(t_1+\pi,t+\pi)$ coincide for $t\in [t_1,t^\ast]$ and differ
in any neighbourhood of $t^\ast.$
 Consequently, there exists a point $(\alpha,\beta)$ on
 $z_K(t^\ast,t_2)$ which is arbitrarily close to $z_K(t^\ast)$ and
 is not in $R z_K(t_1 + \pi, t_2 + \pi).$ We pick the point
 $(\alpha',\beta)$ on $R z_K(t_1+\pi,t_2+\pi),$ which has the same
 ordinate as the point $(\alpha,\beta)$ on $z_K(t_1,t_2)$ (see also Fig.~\ref{06.09.22,10:43}).

In the case $\alpha<\alpha'$ the arc
$z_K(t_1+\pi,t_2+\pi)+(\alpha,\beta)-z_K(t_1+\pi)$ contains the point
$(\alpha,\beta)$ and the point
\begin{equation*} R(\alpha',\beta)+(\alpha,\beta)-z_K(t_1+\pi)
= z_K(t_1)+(\alpha-\alpha',0)= (\alpha-\alpha',0).
\end{equation*}
Since
$(\alpha,\beta) \in z_K(t_1,t_2)$ and $\alpha-\alpha'<0,$ we see that the
endpoint $z_K(t_1)$ of $z_K(t_1,t_2)$ is captured by a translate
of the arc $z_K(t_1+\pi,t_2+\pi).$ The corresponding translation
vector $(\alpha,\beta)-z_K(t_1+\pi)$ can be chosen arbitrarily close to
the vector $z_K(t^\ast)-z_K(t_1+\pi).$ Using the symmetry of the
arcs $z_K(t_1,t^\ast)$ and $z_K(t_1+\pi,t^\ast)$ we get that
$z_K(t^\ast)-z_K(t_1+\pi)=z_K(t_1)-z_K(t^\ast+\pi).$ Hence
$$ z_K(t^\ast)-z_K(t_1+\pi) = \frac{1}{2} (z_K(t^\ast)-z_K(t_1+\pi)+z_K(t_1)-z_K(t^\ast+\pi) ) =\frac{1}{2} ( z_{DK}(t^\ast)+z_{DK}(t_1) ) = z_0.$$

In view of the invariance of the statement of the lemma with
respect to interchanging $t_1$ and $t_2$ with $t_1+\pi$ and
$t_2+\pi,$ respectively, the opposite case $\alpha'<\alpha$ is
settled analogously.
\end{proof}

The following lemma is a strengthening of Proposition~5.1 from \cite{MR1938112}. It states that      $\Ccal_+^2$ regularity of a planar convex body $K$ can be recognized from the covariogram of $K$ over every neighbourhood of $\bd DK.$

\begin{lemma} \label{05.08.10,14:25}
    Let $K$ and $H$ be plane convex bodies with $K \in \Ccal_+^2.$ Then the condition $\GG(\bd DK)$ implies
    that also $H$ belongs to the class $\Ccal_+^2.$
\end{lemma}
\begin{proof}
    Clearly, under the given assumptions we get the equality
$DK=DH.$ Strict convexity of $K$ is equivalent to strict convexity
of $DK.$ Thus, since $K$ is strictly convex, we get that $H$ is
strictly convex, as well. It can be seen that $H$ belongs to
$\Ccal^1$ by examining the asymptotic behaviour of $g_H(x)$
(restricted to $DH$) at boundary points of $H$ (see
\cite[p.190]{MR1938112}).  Further on, in order to get that $H$ is
from the class $\Ccal_+^2$ we can argue in the same way as in the
proof of Lemma~6.2 from \cite{MR2108257}, where the equality of
covariograms is used only at points close to their support.
\end{proof}

\begin{proof}[Proof of Theorem~\ref{05.08.10,15:35}]
Assume that $X=X_0 \cup \bd DK.$ First, by Lemma~\ref{05.08.10,14:25} we deduce that $H \in \Ccal_+^2.$ If $K$ is
centrally symmetric, then the knowledge of the mapping $u \mapsto
\{\tau_K(u),\tau_K(-u)\}, u \in \usphere^1,$ determines $K.$ Thus,
in view of Lemma~\ref{05.08.10,10:47}, we get the assertion.
 Now let us assume that $K$ is not centrally symmetric.  
Further on, let $X'$ be an arbitrary open set with $X \subseteq X'.$ 
Let us prove by contradiction that $K$ is determined within the class $\Ccal_+^2$ by its covariogram over $X'.$ 
Assume the contrary, i.e., there exists a planar convex body $H$ from the class $\Ccal_+^2$ such that $H$ cannot be obtained from $K$ by reflection or translation and $g_K(x)=g_H(x)$ for all $x \in X'.$ 
Let $t_0 \in \real$ be such that $\tau_K(t_0) \ne \tau_K(t_0+\pi).$
In \cite[pp.~186-187]{MR1938112} it is shown that replacing $H$ by an appropriate translation or reflection  there exist arcs $A^+$ and $A^-$ containing $z_K(t_0)$ and $z_K(t_0+\pi),$ respectively, in their relative interiors and contained in the set $\bd K \cap \bd H.$ 
Furthermore, in \cite{MR1938112} it is also noticed that if we assume additionally that $A^+$ and $A^-$ are maximal arcs with the above properties, then $A^+$ and $A^-$ are antipodal to each other,
 i.e., $A^+=z_K(t_1,t_2)$ and $A^-=z_K(t_1+\pi,t_2+\pi)$ for $t_1, t_2 \in \real$ with $t_1< t_0 < t_2$ and $t_2-t_1<\pi.$ 
We shall get a contradiction by showing that $A^+$ or $A^-$ is not maximal, i.e., there exists an arc which strictly contains $A^+$ or $A^-$ and is contained in $(\bd K) \cap (\bd H).$

Obviously, $\tau_K(t)=\tau_H(t)$ for $t \in [t_1,t_2] \cup
[t_1+\pi,t_2+\pi].$ By Lemma~\ref{05.08.10,10:47} we get the
equality
\begin{equation} \label{05.08.19,11:48}
    \{ \tau_K(t), \tau_K(t+\pi) \} = \{ \tau_H(t), \tau_H(t+\pi)\},
\end{equation}
for each $t \in \real.$ If $\tau_K(t_1) \ne \tau_K(t_1+\pi),$ then
using the equalities $\tau_K(t_1)=\tau_H(t_1),$
$\tau_K(t_1+\pi)=\tau_H(t_1+\pi),$ \eqref{05.08.19,11:48} and the
continuity of the functions $\tau_K(t)$ and $\tau_H(t)$ we get
that there exists an $\varepsilon>0$ such that the equality
$\tau_K(t)=\tau_H(t)$ holds for $t \in [t_1-\varepsilon,t_1].$
Consequently, by \eqref{05.02.22,11:26} we have
$z_K(t_1-\varepsilon,t_1)=z_H(t_1-\varepsilon,t_1),$ a
contradiction to the maximality of $A^+.$ Thus, in the sequel we
assume that $\tau_K(t_1)=\tau_K(t_1+\pi).$ If there exists an
$\varepsilon>0$ such that $\tau_K(t)=\tau_K(t+\pi)$ for $t \in
[t_1-\varepsilon, t_1],$ then in view of \eqref{05.08.19,11:48} we
have $\tau_K(t)=\tau_H(t)$ for $t \in [t_1-\varepsilon,t_1],$
which, by \eqref{05.02.22,11:26}, implies the equality of the arcs
$z_K(t_1-\varepsilon,t_1)$ and $z_H(t_1-\varepsilon,t_1),$ a
contradiction to the maximality of $A^+.$

Now let us switch to the case when for every $\varepsilon>0$ the
functions $\tau_K(t)$ and $\tau_K(t+\pi)$ restricted to
$[t_1-\varepsilon,t_1]$ are not identically equal, i.e., there
exists a $t \in [t_1-\varepsilon,t_1]$ with $\tau_K(t) \ne
\tau_K(t+\pi).$ Let $t^\ast$ be the maximal scalar such that $t_1
\le t^\ast \le t_2$ and $\tau_K(t) = \tau_K(t+\pi)$ for $t\in
[t_1,t^\ast].$ If $t^\ast=t_1,$ we put $v=z_K(t_1).$ If
$t^\ast>t_1,$ then for some $n \in \NN$ the arc
$z_K(t_1,t^\ast)$ is a translate of $\frac{1}{2} A_n$ or $- \frac{1}{2}A_n.$ In this case we
put $v=x_n.$

 By Lemma~\ref{05.08.10,14:15} we see that either the endpoint
$z_K(t_1)$ of $A^+$ can be captured by $A^-$ or the endpoint
$z_K(t_1+\pi)$ of $A^-$ can be captured by $A^+.$ Furthermore, the
corresponding translation vector can be chosen arbitrarily
close to $v$ or $-v.$ Without loss of generality, we assume that
$z_K(t_1)$ is captured by $A^-.$ In~\cite[pp.188-189]{MR1938112}
it is shown that in this case a small arc
$z_K(t_1+\pi-\varepsilon', t_1+\pi), \ \varepsilon' >0,$ is
determined by the knowledge of $z_K(t_1,t_2), z_K(t_1+\pi,
t_2+\pi)$ and the values of the covariogram functions at points
arbitrarily close to $z_n.$ This means that we have the equality
$z_K(t_1+\pi-\varepsilon',t_1+\pi)=z_H(t_1+\pi-\varepsilon',t_1+\pi),$
a contradiction to the maximality of $A^-.$

Theorem~\ref{06.05.15,12:36} (Part~III) and the statement of
Theorem~\ref{05.08.10,15:35} for the case $X=X_0\cup \bd DK$
trivially imply the statement of Theorem~\ref{05.08.10,15:35} for
the case $X=X_0 \cup \{o\}.$
\end{proof}

\begin{proof}[Proof of Theorem~\ref{06.05.15,12:36} (Part~II)] We only need to verify the implication $\GG \Rightarrow \LC$ for planar $\Ccal_+^2$ convex bodies $K,$ since the reverse implication is covered by Part~I of the theorem.
We borrow the notations from the statement of Theorem~\ref{05.08.10,15:35}. Let $X'$ be an arbitrary open set with $\bd DK \subseteq X'.$ The set  $X_0$ does not have accumulation points in $\interior DK,$ because for any $\eps>0$ finitely 
many local symmetries of $K$  have length greater than $\eps.$ Therefore only finitely many 
points of $X_0$ lie outside $X'$. Let $\{A^+_i, A^-_i\}$, for $i=1,\dots,n$ and some $n \in \natur,$ be all the 
local symmetries of $K$ corresponding to points of $X_0$ outside $X'$. Let $B^+$ and $B^-$ 
be some antipodal connected components of 
$$\bd K\setminus \bigcup_{i=1}^n(A^+_i\cup A^-_i).$$
Choose $m\in \{1,\ldots,n\}$ such that $A^+_m$ is adjacent to $B^+$ and $A^-_m$ is adjacent 
to $B^-$ (or vice versa).

Let $p$ and $q$ be antipodal points of $\bd K$. Either $p\in \relint A^+_i$ and $q\in 
\relint A^-_i$, for some $i\in\{1,\ldots,n\}$, or  $p\in \relint B^+$ and $q\in \relint 
B^-$, for a suitable choice of $B^+$ and $B^-$, or  $p\in \relint (B^+\cup A^+_m)$ and 
$q\in \relint (B^-\cup A^-_m)$, for a suitable choice of $B^+$ and $B^-$ and $m$. In the 
first case, by Lemma~\ref{05.08.10,10:47}, $\bd K$ and $\bd H$, suitably translated, coincide in a 
neighbourhood of $p$ and $q$. To deal with the second case we observe that the points of 
$X_0$, which correspond to local symmetries contained in $B^+\cup B^-$, belong to $X'$. 
Therefore arguments similar to those of the proof of Theorem~\ref{05.08.10,15:35} (for the case when $X$ from that theorem is given by $X= \bd DK \cup X_0$) imply that if 
$g_K(x)=g_H(x)$ for each $x\in X'$, then $B^+\cup B^-$ is contained in a translate or a 
reflection of $\bd H$. The third case follows from the previous two.
\end{proof}

\begin{remark}\label{06.05.16,16:17}
It is natural to look for minimal (with respect to inclusion) sets $X$ such that $\GG(X)$
implies coincidence of $K$ and $H$, up to translations and reflections. Since the
covariogram is $o$-symmetric, we limit our discussion to $o$-symmetric sets $X.$ We claim that  the set $X$ defined in Theorem~\ref{05.08.10,15:35} is minimal in the following sense. For
certain $\Ccal^2_+$ sets $K$  it suffices to remove from $X_0$ two pairs of opposite points to
violate the conclusion of the theorem. Let us construct a corresponding counterexample. Assume that  two local symmetries of $K$ have
the same center of symmetry, say $o$. Let $\pm A_1$ and $\pm A_2$ be the arcs that
constitute these local symmetries, and $\pm x_1,\pm x_2$ be the midpoints defined 
as in the statement of Theorem~\ref{05.08.10,15:35}, which correspond to the local symmetries $\pm A_1$ and 
$\pm A_2.$ Let $B_i,\ i \in \{1,2,3,4\}$, be the connected components (in
counterclockwise order) of $\bd K\setminus (\pm A_1\cup\pm A_2)$. We
claim that there exist a convex body $H$ which is not a translation or a reflection of
$K$ and such that $\GG(X\setminus\{\pm x_1,\pm x_2\})$ holds. It suffices to define $H$ as
the body obtained from $K$ by flipping the boundary arcs $B_1$ and $B_3.$ That is, the boundary of $H$ is composed of the arcs $\pm A_i \ (i=1,2)$ and $B_2, B_4, -B_1, -B_3.$ The bodies $K$
and $H$ satisfy $\LC$ and thus their covariograms coincide in a neighbourhood of $\bd
DK\cup\{ o \}$, by Theorem~\ref{06.05.15,12:36} (Part~III). Moreover if $C$ and $D$ are the arcs which  constitute a
local symmetry of $K$, different from $\pm A_1$ and $\pm A_2$, then the boundaries of $K$
and $H$ (properly translated and, possibly, reflected) coincide in a neighbourhood of
$C$ and $D$. Therefore $g_K$ and $g_H$ coincide in a neighbourhood of the midpoints
corresponding to the local symmetry.

    We emphasize that the example constructed here is similar in nature to the example from the proof of Theorem~\ref{06.05.08,17:28} (for the case $d=2$).
\end{remark}

    If $K$ is a $\Ccal_+^2$ globally non-symmetric planar convex body, then we see that the set $X_0$ introduced in the proof of Theorem~\ref{05.08.10,15:35} is empty. This remark obviously yields Part~II of Corollary~\ref{06.05.16,17:58}.

\section{Determination results for symmetric and locally \\ symmetric bodies}

\begin{lemma}\label{06.05.08;14:43}
    let $K$ and $H$  be  strictly convex bodies in $\E^2.$ Let $K$ be $o$-symmetric and let $DK=DH.$  Then for scalars $t_0,t_1, t_2 $ with $t_0 \le t_1 \le t_2$ and $t_2 - t_0 \le \pi$ the subset
    $$
        z_H(t_1,t_2) \cup z_H(t_1+\pi,t_2+\pi) \cup \{z_H(t_0),z_H(t_0+\pi)\}
    $$
    of $\bd H$ is centrally symmetric if and only if for every $t \in [t_1,t_2]$ we have $\nabla g_K(x_t)=\nabla g_H(x_t),$ where $x_t:= \frac{1}{2} (z_{DK}(t_0)+z_{DK}(t)).$
\end{lemma}
\begin{proof}
    Let us prove the sufficiency. 
The equality $\nabla g_K(x_t) = \nabla g_H(x_t)$ is equivalent to the condition that $P_H(x_t)$ is a translate of $P_K(x_t)$ (see \eqref{06.09.21,17:07} for the relation among the gradient and $P_K$). 
But, since $K$ is $o$-symmetric, $P_K(x_t)=\conv \{\pm \frac{1}{2} z_{DK}(t_0),\pm \frac{1}{2} z_{DK}(t)\}.$ 
Consequently, the diagonals of $P_H(x_t)$ are translates of  $[o,z_{DK}(t_0)]$ and $[o,z_{DK}(t)].$ 
The chord $[z_H(t_0),z_H(t_0+\pi)]$ is the only chord of $H$ being a translate of $[o,z_H(t_0)-z_H(t_0+\pi)]$ (because that chord is an \emph{affine diameter,} that is, $z_H(t_0)$ and $z_H(t_0+\pi)$ are antipodal). 
Since
    $$
	z_H(t_0)-z_H(t_0+\pi) = z_{DH}(t_0)=z_{DK}(t_0)
    $$
    $[z_H(t_0),z_H(t_0+\pi)]$ is a diagonal of $P_H(x_t).$ 
A similar argument implies that $[z_H(t),z_H(t+\pi)]$ is the other diagonal of $P_H(x_t).$  
Hence, for every $t \in [t_1,t_2],$ the point $z_H(t)$ is a reflection of $z_H(t+\pi)$ with respect to the midpoint $\frac{1}{2} (z_H(t_0)+z_H(t_0+\pi)$ of the diagonal $[z_H(t_0),z_H(t_0+\pi)]$ and the sufficiency is verified.

Now let us show the necessity. 
If $z_H(t_1,t_2)\cup z_H(t_1+\pi,t_2+\pi) \cup \{z_H(t_0),z_H(t_0+\pi)\}$ is centrally symmetric, then a translate of this set is contained in $\frac{1}{2} \bd DH.$ Since $DK=DH$ and $K$ is $o$-symmetric, a translate of the same set coincides with  $z_K(t_1,t_2)\cup z_K(t_1+\pi,t_2+\pi) \cup \{z_K(t_0),z_K(t_0+\pi)\}.$ 
The latter implies that for every $t \in [t_1,t_2]$ the parallelogram $P_H(x_t)$ is a translate of $P_K(x_t)$ and, in view of \eqref{06.09.21,17:07}, shows the sufficiency.
\end{proof}

\begin{proof}[Proof of Theorem~\ref{06.05.16,15:21}]
 \emph{Part~I.} The equality $DK=DH$ implies that $H$ is strictly convex. 
Let $G$ be any open set containing $A.$ We pick an arbitrary $s \in \real$ and show that for a sufficiently small $\eps>0$ the boundary arcs $\bd H \cap (z_{H}(s)+\eps \cdot \uball^2)$ and $\bd H \cap (z_{H}(s+\pi)+\eps \cdot \uball^2)$ 
around the antipodal points $z_{H}(s)$ and $z_{H}(s+\pi),$ respectively, are symmetric with respect to a reflection in a point.
 If $t$ is ranging from $s$ to $s+\pi,$ then the midpoint of the chord $[z_{DH}(s),z_{DH}(t)]$ of $DH$ traverses a path starting at $z_{DK}(s)$ and terminating at the origin. 
Thus, for some $t_0 \in [s,s+\pi]$ the midpoint $\frac{1}{2} (z_{DH}(s)+z_{DH}(t_0))$ of $[z_{DH}(s),z_{DH}(t_0)]$ lies in $A.$

         \begin{figtabular}{c}
        \begin{picture}(40,30)
            \put(4,4){\ExternalFigure{0.30}{057/057n34.eps}} 
        \put(36,14){$DH$}
        \put(21,16){$o$}
	\put(9,28){$z_{DH}(s)$}
	\put(-1,8){$z_{DH}(t_0)$}
        \end{picture}
        \\
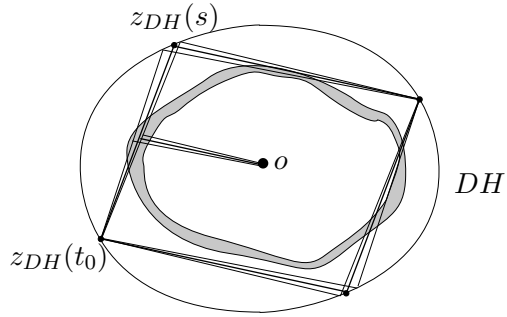
    \parbox[t]{0.60\textwidth}{\caption{The set $G$ containing $A$ is painted in gray}}
    \end{figtabular}

 Clearly, for some $\eps>0$ the midpoint of $[z,z_{DH}(t_0)]$ lies in $G$ for all $z \in \bd DH$ with $|z-z_{DH}(s)|<\eps.$ Let $t_1,t_2$ be scalars such that
 $$
    \bd DH \cap (z_{DH}(s)+\eps \cdot \uball^2) = \setcond{z_{DH}(t)}{t_1 \le t \le t_2}.
 $$
 The assumption of the theorem and \eqref{06.09.21,17:07} imply that $\nabla g_H(x_t) = \nabla g_K(x_t),$ where $x_t := \frac{1}{2} (z_{DH}(t_0)+z_{DH}(t))$ and $t \in [t_1,t_2].$ 
Hence, in view of Lemma~\ref{06.05.08;14:43}, $z_H(t_1,t_2)$ is a reflection of  $z_K(t_1+\pi,t_2+\pi)$ in a point. 
Thus, each pair of antipodal points of $H$ can be enclosed in the relative interior of symmetric boundary arcs of $H,$ which implies the central symmetry of $H.$

 \emph{Part~II.} The case $X=\{o\}$ can be transformed to the case $X=\bd DK$ using Part~III of Theorem~\ref{06.05.15,12:36}. 
Thus, we assume that $X=\bd DK.$ The statement of the theorem follows then from Lemma~\ref{06.05.08;14:43} applied for arbitrary $t_0 \le t_1 \le t_2$ with $t_0=t_1$ and $t_2$ sufficiently close to $t_1.$

 \end{proof}

 \begin{proof}[Proof of Corollary~\ref{06.05.24,09:59}]
 Let $K$ be an arbitrary locally symmetric convex body in $\E^2.$ Let us show that $g_K$ determines $K,$ up to translations and reflections. 
We assume that $K$ is strictly convex, since for non-strictly convex bodies the determination was verified in \cite[Theorem~1.1]{MR2108257}. 
By Part~II of Theorem~\ref{06.05.16,15:21} we can determine all local symmetries of $K,$ up to translations. 
But then the theorem follows from Proposition~1.4 in \cite{MR2108257}, stating that if additionally to the knowledge of $g_K$ a non-degenerate boundary arc of $K$ is known, then $K$ can be determined uniquely, up to translations and reflections.
 \end{proof}

\section{Genericity results}

The \emph{Nikodym distance} $\delta_N$ (also known as the 
\emph{symmetric difference metric})  between convex bodies 
$K$ and $H$ in $\E^d$ is given by
\begin{equation}
 \delta_N(K,H)= V ( (K\setminus H) \cup (H\setminus K)) .
\end{equation}
It is known that $\delta_N$ generates the same topology in the
class of convex bodies as the Hausdorff distance $\delta,$
\cite[pp.58-59]{MR94d:52007}. Furthermore
 the inequality (see \cite[p.195]{MR1938112})
\begin{equation} \label{05.09.01,11:40}
    |g_K(x)-g_H(x)| \le 2 \delta_N(K,H)
\end{equation}
for all $x\in\E^d$ shows that the operator $K \mapsto g_K$ is continuous provided the class of
convex bodies is endowed with the Nikodym distance, and the
distance between covariograms is measured with respect to the
maximum norm.

\begin{lemma} \label{05.08.12,11:37}
The class of totally non-symmetric $\Ccal^2_+$ planar convex bodies is dense in the class
of all  $\Ccal^2_+$ planar convex bodies, with respect to the Hausdorff metric.
\end{lemma}
\begin{proof}
Let $R(u):=R_K(u).$
It suffices to approximate a $\Ccal^2_+$ convex body $K$ by a totally non-symmetric $\Ccal^2_+$ one. In order to do this we approximate,  for each $\varepsilon>0$,  $R(u)$ by a continuous positive function $R_\varepsilon(u)$ with the property that the set $U:=\setcond{u\in \usphere^1}{R_\varepsilon(u)=R_\varepsilon(-u)}$ has empty relative interior and
\begin{align*}
    &\int_{\usphere^1} R_\varepsilon(u) \cdot u  \, \dd u = 0,&
&\max_{u\in \usphere^1}|R(u)-R_\varepsilon(u)|<\varepsilon.
\end{align*}
By Minkowski's theorem and the relative  stability result (Theorem~7.2.2 from \cite{MR94d:52007}) there exists a convex body $K_\varepsilon,$ whose radius of curvature is $R_\varepsilon(u)$ and
$\bar{\delta}(K,K_\varepsilon) = O(\varepsilon).$

Let us construct $R_\eps$. The
relative interior of $U$, $\relint U$, is open and therefore it is the
disjoint union of denumerably many open intervals. It is thus possible
to construct an odd continuous function $f(u)$ in $\usphere^1$ with the
property that $|f(u)|<1$ for each $u$,  $f(u)$ vanishes outside $\relint
U$, and  $f(u)$ never vanishes in $\relint U$ (except for the case
$U=\usphere^1$, since in this case  it has to vanish in at least two antipodal
points). The function $R_{\eps}=R+ \eps f$, for $0<\eps<\min_{u\in \usphere^1}
R(u)$, satisfies the required properties.
\end{proof}

Now we are ready to give the proof of the genericity statement
given in Theorem~\ref{06.05.08,17:28} (Part~I). We shall settle the cases
$d=2$ and $d \ge 3$ independently of each other.

\begin{proof}[Proof of Theorem~\ref{06.05.08,17:28} (Part~I) for $d=2$]
Let us denote by $\KK'$ the class of all planar convex bodies $K$ which are not determined, up to translations and reflections, by their covariogram over every neighbourhood of $\bd DK.$
Given $n \in \natur$ we introduce the class $\KK_n' \subseteq \KK^2_0$ such that $K \in \KK_n'$ if and only if there exists a planar convex body $H$ with the properties
\begin{align*}
    &DK = DH,&
    &g_K(x) = g_H(x) \quad \text{for $x\in\bd DK+\frac{1}{n} \cdot \uball^2$,} \\
    &\bar{\delta}(K,H) \ge \frac{1}{n},&
    &H \subseteq  n \cdot \interior \uball^2
\end{align*}
Clearly $\KK' = \bigcup_{n=1}^{+\infty} \KK_n'$.
Let us prove now that for every $n \in \natur$ the class $\KK_n'$ is closed. Let $(K_m)_{m=1}^{+\infty}$ be a sequence of convex bodies belonging to $\KK_n'$ and converging to $K_0,$ and for each $m$
 let $H_m$ be the  convex body associated to $K_m$.
By the Blaschke selection theorem there exists a subsequence $(H_{m_j})_{j=1}^{+\infty}$ converging to some $H_0 \in \KK^2.$
 Obviously, $H_0 \subseteq n \cdot \interior\uball^2$, $\bar{\delta}(K_0,H_0) \ge \frac{1}{n}$ and $D K_0= D H_0.$

Let $x$ be a point from $\bd DK_0 + \frac{1}{n} \cdot \interior \uball^2.$
Then $x$ belongs to  $\bd DK_{m_j} + \frac{1}{n} \cdot \interior \uball^2$
definitely, i.e.\  for all $j$ sufficiently large. Since $g_{K_{m_j}}(x) = g_{H_{m_j}}(x)$ and passing to the limit (recall \eqref{05.09.01,11:40})  one obtains $g_{K_0}(x) = g_{H_0}(x)$. Therefore $K_0$ is indeed in $\KK_n'.$

Thus, $\KK'$ is a countable union of closed sets. Consequently,
the complement $\Ucal$ of $\KK'$ is a countable intersection of
open sets. Using Corollary~\ref{06.05.16,17:58} (Part~I), we get
that $\Ucal$ contains all $C_+^2$ totally non-symmetric planar
convex bodies. Further on, applying Lemma~\ref{05.08.12,11:37},
wee see that $\Ucal$ is dense in the class of all convex bodies in
$\E^2.$
\end{proof}

\newcommand{\HH}{\mathcal{H}}
\newcommand{\Sc}{\mathcal{S}}
\begin{proof}[Proof of Theorem~\ref{06.05.08,17:28} (Part~II)]
Let $\HH'$ be the class of all planar convex bodies $K$ which are not determined, up to translations and reflections, by their
    covariogram over every neighbourhood of the origin. 
Let $\KK'$ be as in the previous proof, and $\Sc$ be the class of all planar strictly convex bodies. 
We have $\HH' = (\HH' \setminus \Sc) \cup (\HH' \cap \Sc).$  It is well known that the class $\KK_0^2 \setminus \Sc$ of all non-strictly convex bodies is meager. 
Therefore, its subclass $\HH' \setminus \Sc$ is meager, as well. The class $\HH'\cap \Sc$ is meager because, by Theorem~\ref{06.05.15,12:36} (Part~III) , it is a subclass of
     $\KK'$. 
Then $\HH'$ is meager, since it is the union of two meager classes.
\end{proof}

Given convex bodies $K, H \subseteq \E^d,$ the \emph{cross covariogram function} of $K$ and $H$ is defined by
\begin{equation*}
    g_{K,H}(x) := V(K \cap (H+x)),
\end{equation*}
where $x$ ranges over $\E^d.$ The support of $g_{K,H}(x)$ is
obviously equal to $K+(-H).$ For the proof of Theorem~\ref{06.05.08,17:28} (Part~I) for $d \ge 3$
we need the following lemma.

\begin{lemma} \label{05.10.24,13:17}
    Let $T$ be a $d$-dimensional simplex in $\E^d.$ Let $\lambda \in (0,1).$ Then the following statements hold true.
    \begin{enumerate}[I.]
        \item If $p$ is a vertex of $T$ and $-u \in S^{d-1}$  is in the relative interior of the support cone of $T$ at $p,$ then for sufficiently small $t>0$
        \begin{equation} \label{06.09.21,15:33}
            g_{(1-\lambda)T,-\lambda T}(p-t u) = C t^d,
        \end{equation}
        where $C=C(T,p,u)$, i.e. $C$  does not depend on $\lambda$.
        \item If $I$ is an edge of $T,$ $p$ is the midpoint of $I,$ and $-u \in S^{d-1}$ is in the relative interior of the support cone of $T$ at $p,$ then for sufficiently small $t >0$
        \begin{equation} \label{06.09.21,21:32}
            g_{(1-\lambda)T,-\lambda T}(p-t u) = C \min \{1-\lambda,\lambda\} t^{d-1} + o(t^{d-1}),
        \end{equation}
        where $C=C(T,I,u).$
    \end{enumerate}
\end{lemma}
\begin{proof}
    \emph{I.} Let $P_t:=[(1-\lambda)] T \cap [-\lambda T+p-tu].$ 
Since $(1-\lambda)p$ is a vertex of $(1-\lambda)P,$ it is also a vertex of $P_t$ provided $(1-\lambda) p  \in -\lambda T + p- t u.$ The latter condition is equivalent to $p-\frac{t}{\lambda} u \in T.$ 
Analogously, the vertex $(1-\lambda) p - tu$ of $-\lambda T + p -tu$ is also a vertex of $P_t$ provided $p- \frac{t}{1-\lambda} u \in T.$ 
Let $t_0>0$ be such that 
point $p-1/\min\{1-\lambda,\lambda\} \cdot u$ lies in $T$. Then in view of the above remarks for $0 < t < t_0$ both $(1-\lambda)p$ and $(1-\lambda) p - tu$ are vertices of $P_t.$ 
Moreover, it can be seen that $P$ is a parallelotope whose facets are parallel to facets of $T$ incident to $p,$ and $[(1-\lambda) p, (1-\lambda)p - tu]$ is a diagonal of $P_t.$ 
Then $t u$ is a vector joining the endpoints of this diagonal. 
Consequently,
for $t_1, t_2 \in [0,t_0]$ the polytopes $P_{t_1}$ and $P_{t_2}$
are homothetic with homothety ratio $\frac{t_1}{t_2}$ and $P_t$
does not depend on the choice of $\lambda.$ The above facts easily
imply the statement of Part~I.

    \emph{II.} Let $P_t$ be introduced in the same way as above, and $L$ be a hyperplane through the origin orthogonal to  the line $\aff I.$  
Then the set $P_t$ can be approximated by the cylinder  $P_t \ortproj L + \min \{1-\lambda,  \lambda \} I $ so that we have
    \begin{align*}
    g_{(1-\lambda)T,-\lambda T }(p-tu) &=  V( (P_t \ortproj L) + \min \{1-\lambda,\lambda\}  I)+ o(t^{d-1})  \\ &= V_{d-1} (P_t \ortproj L) \cdot
        \min \{1-\lambda,\lambda\} \cdot V_1(I) + o(t^{d-1}).
    \end{align*}
Clearly, for small $t$ the polytope $P_t \ortproj L$ is the intersection of the $(d-1)$-dimensional simplices $(1-\lambda) T  \ortproj L$ with $(-\lambda  T + p-tu) \ortproj L,$ which means
 that we can apply the statement of Part~I for the cross covariogram function of these $(d-1)$-dimensional simplices. 
By this we get the statement of Part~II.
\end{proof}

\begin{proof}[Proof of Theorem~\ref{06.05.08,17:28} (Part~I) for $d \ge 3$]
 As was mentioned in \cite{MR1416411}, it is sufficient to prove
the determination property for  $g_P(x)$ in the case when $P$ is a
simplicial polytope such that  $P$ and $-P$ are in
general relative position. It is known (see
\cite[p.86]{MR1416411} and \cite[Theorem~2.1]{MR1322067}) that if
a polytope $H$ has the difference body $DP$ (for $P$ as above),
then $H=(1-\lambda) P + \lambda (-P).$ Thus, clearly for the
 determination of $H,$ up to translations and reflections, it is sufficient to retrieve the set
$\{1-\lambda, \lambda\}.$  Let $u$ be an outward facet normal of
$P$ and let $T:=F_P(u).$ Since $P$ is a simplicial polytope, $T$
is a $(d-1)$-dimensional simplex. Since $P$ and $-P$ are in
general relative position we have that $F_P(-u)$ is a singleton
consisting of some vertex $p$ of $P.$  Clearly, $(1-\lambda) p -
\lambda T=F_H(-u)$ and $(1-\lambda)T-\lambda p=F_H(u)$. Let us
pick a point $x^+$ from $F_H(u)$ and a point $x^-$ from $F_H(-u).$
Consider the point $x:=x^+ - x^-$ from $T-p$ (which is the facet
of $DP$ parallel to $T$). Let us consider an arbitrary vector
$y \in u^\perp.$ Then, for small $t>0$
\begin{equation} \label{06.09.21,21:56}
    g_H(x-tu+y) = V( H \cap (H+x-tu+y)  ) =  V( (H - x^+) \cap (H-x^- -tu+y)  )
\end{equation}
Polytopes $H-x^+$ and $H-x^-$ involved in \eqref{06.09.21,21:56} have facets $F_H(u)-x^+$ and $F_H(u)-x^-,$ respectively, both lying in $L.$ Thus, it can be seen that the polytope
$$
(H - x^+) \cap (H-x^- -tu+y)
$$
from \eqref{06.09.21,21:56} can be approximated by the cylinder $T_1 \cap (T_2+y)+ [o,-tu],$
where
\begin{align*}
    &T_1 =  (1-\lambda) T - \lambda p - x^+,&
    &T_2 = -\lambda T + (1-\lambda)p - x^-,
\end{align*}
and
\begin{equation}\label{asympt}
   g_H(x-tu+y) = t \cdot g_{T_1,T_2}(y)+o(t).
\end{equation}
Clearly $\lambda \in \{0,1\}$ if and only if either $F_H(u)$ or
$F_H(-u)$ is a point. This is equivalent to $g_{T_1,T_2}\equiv0$,
and in view of \eqref{asympt}, to $g_H(x-tu+y)=o(t)$ for each $y$.
When $\lambda \in (0,1)$ then the function $g_{T_1,T_2}$ is
determined by \eqref{asympt}. In view of
Lemma~\ref{05.10.24,13:17} (Part~II), this function determines
the set $\{1-\lambda,\lambda\}$.
\end{proof}

Given a convex body $K \in \KK_0^d$ and a vector $x \in \E^d,$ we
introduce the body
\begin{equation*}
    K(x):= \setcond{y \in K}{ V_1(K \cap \aff \{y,y+x\}) \ge |x|},
\end{equation*}
which is the union of all those chords of $K$ that are parallel to
$x$ and are not shorter than $[o,x].$ It can easily be shown that
$K \cap (K+x) = K(x) \cap (K(x) + x).$ Consequently,
\begin{equation} \label{05.09.25,20:15}
          g_{K(x)}(x)=g_K(x).
\end{equation}

Let $$\bar{K}(x):= \setcond{y \in x^\perp}{V_1(K\cap \aff \{y,y+x\}
) \ge |x|}.$$ Clearly, $\bar{K}(x)$ is the orthogonal
projection of $K(x)$ onto $x^\perp.$ For $u \in \usphere^{d-1}$ we
have
\begin{equation} \label{05.12.27,14:40}
    \frac{\partial}{\partial t} g_K(t u) = -V_{d-1}(\bar{K}(t u)).
\end{equation}
Formula \eqref{05.12.27,14:40} is  presented in \cite{Matheron8601} and
\cite{MR0385969}.
Now let us come to the proof of the next theorem.
\begin{proof}[Proof of Theorem~\ref{06.05.08,17:28} (counterexample)]
The counterexample constructed below is strongly related to some counterexample from \cite{MR1416411}.
 Let $U_1$ and $U_2$ be relatively open subsets of $\usphere^{d-1}$ bounded by
$(d-2)$-dimensional spheres and such that the sets $\pm U_1, \pm
U_2$ are mutually disjoint. Let $K$ be a $\Ccal_+^2$ convex body
satisfying the conditions $r_K(u)=1$ for $u \in \usphere^{d-1}
\setminus ( U_1 \cup U_2)$ and $r_K(u)<1$ for $u \in (U_1\cup
U_2).$ Then we introduce the body $H$ defined by $r_H(u):=r_K(-u)$
for $u \in (-U_1) \cup U_1 $ and $r_H(u):=r_K(u),$ otherwise (see
Figs.~\ref{05.08.30,16:50}-\ref{05.08.30,16:51} for the
illustration in the case $d=2$). It can be seen that $DK=DH.$

Let
\begin{align*}
    &A_1 :=  \setcond{z_K(u)}{u \in \usphere^{d-1}\setminus (U_1 \cup (-U_1))},&
    &A_2 :=  \setcond{z_K(u)}{u \in \usphere^{d-1}\setminus (U_2 \cup (-U_2))}.
\end{align*}
Since $\bd K = A_1 \cup A_2$ and $\bd H = (-A_1) \cup A_2,$  the bodies $K$ and $H$ satisfy $\LC.$ Thus, by
Theorem~\ref{06.05.15,12:36} (Parts~I and II), $g_K=g_H$ in a neighbourhood of their support.

        \begin{figtabular}{ccc}
        \begin{picture}(30,25)
            \put(2,2){\ExternalFigure{0.20}{056/fig056n01v01.eps}} 
            \put(2,22){$K$}
        \end{picture}
        &
        \begin{picture}(25,25)
            \put(2,4.5){\ExternalFigure{0.15}{057/057n32.eps}} 
            \put(2,21.5){$\usphere^1$}
        \end{picture}
        &
        \begin{picture}(30,25)
            \put(2,2){\ExternalFigure{0.20}{056/fig056n01v02.eps}} 
            \put(2,22){$H$}
        \end{picture}
        \\
    \parbox[t]{0.30\textwidth}{\mycaption{\label{05.08.30,16:50}}}
        &
    \parbox[t]{0.25\textwidth}{\mycaption{\label{06.02.23,17:06}}}
        &
    \parbox[t]{0.30\textwidth}{\mycaption{\label{05.08.30,16:51}}}
    \end{figtabular}
 By standard compactness arguments, there exists an $\alpha>0$ such that for every $x \in \E^d$ with $|x| < \alpha$ and for every two-dimensional linear space $L$ containing $x$ the endpoints 
of the two chords of $K\cap L$ which are  translates of $[o,x]$ either all belong to $A_1$ or  all belong to $A_2$.  
Hence, for every $u,v \in \usphere^{d-1}$, with $\sprod{v}{u}=0$ and for $|t|<\alpha$  we have
 \begin{equation} \label{06.02.27,12:47}
        \{r_{\bar{K}(tu)}(v),r_{\bar{K}(tu)}(-v)\} = \{r_{\bar{H}(tu)}(v),r_{\bar{H}(tu)}(-v)\}.
 \end{equation}

We recall that the  volume of a convex body $K$ in $\E^d$ with
$o \in K$ can be written as
\begin{equation} \label{05.12.27,13:23}
    V(K) = \frac{1}{d} \int_{\usphere^{d-1}} r_K(u)^d \,
    \dd u.
\end{equation}
 Let $u \in \usphere^{d-1}$ and $t$ be such that  $0 <t  <\alpha$ and $g_K(tu) >0$. Up to translations of $K,$ we may assume that $o\in \bar{K}(tu),$ and we have
   \begin{equation*} 
\begin{split}
    \frac{\partial}{\partial t} g_K(t u)
    \equalby{\eqref{05.12.27,14:40}} -V_{d-1} ( \bar{K}(t u))
     &\equalby{\eqref{05.12.27,13:23}} \frac{1}{d-1} \int_{\usphere^{d-1} \cap
    u^\perp}r_{\bar{K}(t u)}(v)^{d-1} \, \dd v \\  &= - \frac{1}{2(d-1)} \int_{\usphere^{d-1} \cap
    u^\perp} \left(r_{\bar{K}(t u)}(v)^{d-1} + r_{\bar{K}(t u)}(-v)^{d-1} \right)\, \dd v\\
&\equalby{\eqref{06.02.27,12:47}} \frac{\partial}{\partial t} g_H(t u)
\end{split}
    \end{equation*}
In view of the
 equalities, $V(K)=g_K(o)=g_H(o)=V(H),$ the latter implying the
 coincidence of $g_K$ and $g_H$ for $x \in \E^d$ in a neighbourhood of $o$.
\end{proof}

\def\cprime{$'$} \def\cprime{$'$} \def\cprime{$'$} \def\cprime{$'$}
\providecommand{\bysame}{\leavevmode\hbox to3em{\hrulefill}\thinspace}
\providecommand{\MR}{\relax\ifhmode\unskip\space\fi MR }
\providecommand{\MRhref}[2]{%
  \href{http://www.ams.org/mathscinet-getitem?mr=#1}{#2}
}
\providecommand{\href}[2]{#2}

\begin{tabular}{lr}
\begin{tabular}{l}
  Gennadiy Averkov \\
  Faculty of Mathematics \\
  University of Magdeburg, Universit\"atsplatz 2 \\
  D-39106 Magdeburg \\
  Germany \\
  e-mail: gennadiy.averkov@googlemail.com
\end{tabular}
&

\begin{tabular}{l}
  Gabriele Bianchi \\
  Department of Mathematics \\
  Universit\`{a} di Firenze, Viale Morgagni  67a \\
  50134 Firenze \\
  Italy \\
  e-mail: gabriele.bianchi@unifi.it
\end{tabular}
\end{tabular}

\clearpage
\pagenumbering{roman}

\end{document}